\documentclass[preprint]{elsarticle}
\usepackage{amsmath}
\usepackage{amssymb}
\usepackage{amsfonts}
\usepackage{graphicx}
\usepackage{lineno,hyperref}
\usepackage{subfigure}
\newtheorem{thm}{Theorem}

\newtheorem{proposition}[thm]{Proposition}
\newtheorem{corollary}[thm]{Corollary}
\newtheorem{example}[thm]{Example}
\newdefinition{rmk}{Remark}
\newdefinition{definition}{Definition}
\newproof{pf}{Proof}

\modulolinenumbers[5]

\journal{Journal of Geometry and Physics}

\begin{document}

\begin{frontmatter}

\title{On geometrical structures, associated with linear differential
operators of the 1st order}

\author[mymainaddress]{Valentin Lychagin\corref{mycorrespondingauthor}}
\cortext[mycorrespondingauthor]{Corresponding author}
\ead{valentin.lychagin@uit.no}

\address[mymainaddress]{V.A. Trapeznikov Institute of Control Sciences, Russian Academy of Sciences, 65 Profsoyuznaya Str., 117997 Moscow, Russia}

\begin{abstract}
The problem of equivalency for linear differential operators of the first
order is discussed.
\end{abstract}


\end{frontmatter}

\section{Introduction}

In this paper we continue to study invariants and structures on linear
differential operators. The cases of differential operators of order $k\geq
2,$ acting in line bundles,  were studied in papers (\cite{LY2},\cite{LY3},%
\cite{LYk}).

In this paper we study linear differential operators of the 1st order,
acting in sections of a vector bundle $\pi :E\left( \pi \right) \rightarrow
M,$ where $\dim M=n\geq 2,\dim \pi =m\geq 2.$

To illustrate our approach we consider a "toy case" of the 1-st order
operators , when $\dim \pi =1,$ $\dim M=1.$

Let $\Delta $ be such an operator. Then its symbol $\sigma \in \mathrm{End}%
\left( \pi \right) \otimes \Sigma _{1}=\Sigma _{1}$ is a vector field on $M$.

Assume that we have a connection $\nabla $ in the bundle $\pi $ and let $%
Q_{\nabla }$ be the quantization, associated with $\nabla $ (see, for
example,\cite{LYk} ), i.e., in our case $Q_{\nabla }:\mathrm{End}\left( \pi
\right) \otimes \Sigma _{1}\rightarrow \mathbf{Dif}_{1}\left( \pi ,\pi
\right) $ is the morphism splitting exact sequence $\mathbf{0\rightarrow }%
\mathrm{End}\left( \pi \right) \rightarrow \mathbf{Dif}_{1}\left( \pi ,\pi
\right) \overset{\text{smbl}}{\rightarrow }\mathrm{End}\left( \pi \right)
\otimes \Sigma _{1}\rightarrow \mathbf{0.}$

Then we get a decomposition $\Delta =Q_{\nabla }\left( \sigma \right)
+\sigma _{0}\left( \nabla \right) ,$ where operator $\sigma _{0}\left(
\nabla \right) \in \mathrm{End}\left( \pi \right) $ depends on connection
in the following way.

If $\widetilde{\nabla }=\nabla +\alpha $ is another connection in the bundle
$\pi ,$ where $\alpha $ is a differential $1$-form on $M,$ then $\sigma
_{0}\left( \widetilde{\nabla }\right) =\sigma _{0}\left( \nabla \right)
+\left\langle \alpha ,\sigma \right\rangle .$

Therefore$,\ $if$\ \sigma \neq 0$, we can choose a unique connection $\nabla
$ in such a way, that $\sigma _{0}\left( \nabla \right) =0.$ For this
connection, we have $\Delta =Q_{\nabla }\left( \sigma \right) $ and
equivalence of operators with respect to group of automorphisms of $\pi $ is
equivalent to joint classification of pairs (vector field, connection).

In the case, considered in this paper, we start with classification of
symbols and find a class of symbols, which we call symbols in general
position, such that they allow us to associate with operator a unique
connection (the connection that "preserves"  the symbol in the best way) and
therefore transform the equivalence problem of operators to the equivalence
problem of tensors and connections.

To find invariants of symbols we use the theorem of Procesi (\cite{Pr},\cite%
{Pr2}), where he proved the M.Artin conjecture (\cite{Art}) on\ algebra of
polynomial invariants of a set of linear operators. Applying the Artin's
idea we construct tensor invariants of symbols (we call them Artin-Procesi
invariants) and then describe the field of rational invariants of symbols.
Then, we use these results and ideas of (\cite{LYk}) on natural coordinates
and natural atlases in order to get a local as well as a global description
of regular differential operators of the first order.

The paper is organized in the following way. At first, we discuss
Artin-Procesi invariants of symbols and corresponding conditions of
regularity or general position. We show that the regularity conditions of
symbols provide us with a series of geometrical structures on the manifold.
Namely, they give us a (pseudo) Riemannian structure on $M$ as well as a
canonical frame. Secondly, we study connections in the bundle $\pi $ and on
the manifold, associated with symbols. On the manifold we take the
Levi-Civita connection, associated with the (pseudo) Riemannian structure.
The requirement on connections in $\pi $ to preserve the symbol is too
strong (cf. \cite{LYk}) and we restrict ourselves by the class of
connections (we call them minimal) that "preserve" symbols in the best way.
Using the subsymbol of the operator (similar to  \cite{LYk}) we show that
among minimal there is an unique connection associated with the operator.
Existence of such connection allows us to get local and global description
of operators by using natural coordinates and natural atlases.

\section{Differential operators}

\subsection{Quantizations}

The notations we've used in this paper are similar to notations that were
used in paper (\cite{LYk}) .

Let $M$ be an $n$-dimensional manifold and let $\pi :E\left( \pi \right)
\rightarrow M$ be a vector bundle.

We denote by $\tau :TM\rightarrow M$ and $\tau ^{\ast }:T^{\ast
}M\rightarrow M$ \ the tangent and respectively cotangent bundles over
manifold $M,$ and by $\mathbf{1:}\mathbb{R}\times M\rightarrow M\ $we denote
the trivial linear bundle $.$

The symmetric and exterior powers of a vector bundle $\pi :E\left( \pi
\right) \rightarrow M$ \ will be denoted by $\mathbf{S}^{k}\left( \pi
\right) $ and $\mathbf{\Lambda }^{k}\left( \pi \right) .$

The module of smooth sections of bundle $\pi $ we denote by $C^{\infty
}\left( \pi \right) ,$ and for the cases tangent, cotangent and the trivial
bundles we'll use the following notations: $\Sigma _{k}\left( M\right)
=C^{\infty }\left( \mathbf{S}^{k}\left( \tau \right) \right) -$ the module
of symmetric $k$-vectors and $\Sigma ^{k}\left( M\right) =C^{\infty }\left(
\mathbf{S}^{k}\left( \tau ^{\ast }\right) \right) -$ the module of symmetric
$k$-forms, $\Omega _{k}\left( M\right) =C^{\infty }\left( \mathbf{\Lambda }%
^{k}\left( \tau \right) \right) -$ the module of skew-symmetric $k$-vectors
and $\Omega ^{k}\left( M\right) =C^{\infty }\left( \mathbf{\Lambda }%
^{k}\left( \tau ^{\ast }\right) \right) $ -the module of exterior $k$-forms,
$C^{\infty }\left( \mathbf{1}\right) =C^{\infty }\left( M\right) .$

Let $\mathbf{Diff}_{k}\mathbf{(\pi ,\pi )}$ be the module of linear
differential operators of order $k,$ acting in the sections of the bundle $%
\pi .$ These modules are connected by exact sequences
\begin{equation*}
\mathbf{0\rightarrow Diff}_{k-1}\left( \pi ,\pi \right) \mathbf{\rightarrow
Diff}_{k}\left( \pi ,\pi \right) \overset{\mathrm{smbl}_{k}}{\rightarrow }%
\mathrm{End}\left( \pi \right) \otimes \Sigma _{k}\left( M\right)
\rightarrow \mathbf{0,}
\end{equation*}%
where $\mathrm{End}\left( \pi \right) =\mathbf{Diff}_{0}\left( \pi ,\pi
\right) $ is the module of endomorphisms of $C^{\infty }\left( \pi \right) ,$
and mapping $\mathrm{smbl}$ assign to differential operator $\Delta $ its
symbol $\mathrm{smbl}\left( \Delta \right) .$

If we consider $\mathrm{smbl}_{k}\left( \Delta \right) $ as a linear
operator $\mathrm{smbl}_{k}\left( \Delta \right) :\Sigma ^{k}\left(
M\right) \rightarrow \mathrm{End}\left( \pi \right) ,$ then
\begin{equation*}
\mathrm{smbl}_{k}\left( \Delta \right) \left( df^{k}\right) =\frac{%
1}{k!}\delta _{f}^{k}\left( \Delta \right) \in \mathrm{End}\left( \pi
\right) ,
\end{equation*}%
where $f\in $ $C^{\infty }\left( M\right) ,\ \delta _{f}:\mathbf{Diff}%
_{i}\left( \pi ,\pi \right) \mathbf{\rightarrow Diff}_{i-1}\left( \pi ,\pi
\right) ,$ is the commutator $\delta _{f}\left( \square \right) =f\circ
\square -\square \circ f,$ and $df^{k}=df\cdot \cdots \cdot df\in \Sigma
^{k}\left( M\right) $ is the $k$-th symmetric power of the differential $1$%
-form $df\in \Sigma ^{1}\left( M\right) .$

For the case $k=1,$ which we'll consider in this paper, we respectively have:

\begin{enumerate}
\item Exact sequences%
\begin{equation}
\mathbf{0\rightarrow }\mathrm{End}\mathbf{\left( \pi \right) \rightarrow
Diff}_{1}\left( \pi ,\pi \right) \overset{\mathrm{smbl}_{1}}{\rightarrow }%
\mathrm{End}\left( \pi \right) \otimes \Sigma _{1}\left( M\right)
\rightarrow \mathbf{0,}  \label{Exact1}
\end{equation}%
\newline
\newline
and $\mathrm{smbl}_{1}\left( \Delta \right) \left( df\right)
=\delta _{f}\left( \Delta \right) \in \mathrm{End}\left( \pi \right) .$

\item Any connection $\nabla $ in the vector bundle $\pi $ gives us \textit{%
quantization}
\begin{eqnarray}
\mathbf{Q}_{\nabla } &\mathbf{:}&\mathrm{End}\left( \pi \right) \otimes
\Sigma _{1}\left( M\right) \rightarrow \mathbf{Diff}_{1}\left( \pi ,\pi
\right) ,  \label{ConQuan} \\
&&\mathrm{smbl}_{1}\circ \mathbf{Q}_{\nabla }\mathbf{=}\mathrm{id}%
,  \notag
\end{eqnarray}%
where%
\begin{equation}
\mathbf{Q}_{\nabla }\left( A\otimes X\right) =A\circ \nabla _{X},
\label{Quant}
\end{equation}%
for all $A\in \mathrm{End}\left( \pi \right) ,\ X\in \Sigma _{1}.$
\end{enumerate}

Moreover, if $\ \widetilde{\nabla }$ is another connection and
\begin{equation*}
\widetilde{\nabla }=\nabla +\alpha ,
\end{equation*}%
where $\alpha \in \mathrm{End}\left( \pi \right) \otimes \Sigma ^{1}\left(
M\right) ,$ $\widetilde{\nabla }_{X}=\nabla _{X}+\alpha \left( X\right) ,$
then
\begin{equation*}
\mathbf{Q}_{\widetilde{\nabla }}=\mathbf{Q}_{\nabla }+\left\langle \cdot
,\alpha \right\rangle ,
\end{equation*}%
where $\left\langle ,\right\rangle :\left( \mathrm{End}\left( \pi \right)
\otimes \Sigma _{1}\left( M\right) \right) \otimes \left( \mathrm{End}%
\left( \pi \right) \otimes \Sigma ^{1}\left( M\right) \right) \rightarrow
\mathrm{End}\left( \pi \right) $ is the natural pairing: $\left\langle
A\otimes X,B\otimes \omega \right\rangle =\left\langle \omega
,X\right\rangle \ A\circ B.$

In what follows, we denote by $\sigma _{1}\left( \Delta \right) $ and $%
\sigma _{0}^{\nabla }\left( \Delta \right) \ $(or simply by $\sigma _{1}$
and $\sigma _{0}$) the symbol $\sigma _{1}\left( \Delta \right) =\mathrm{%
smbl}_{1}\left( \Delta \right) $ and \textit{subsymbol} $\sigma
_{0}^{\nabla }(\Delta )=\Delta -\mathbf{Q}_{\nabla }\left( \Delta \right)
\in \mathrm{End}\left( \pi \right) .$

Then symbol $\sigma _{1}$ does not depend on connection but
\begin{equation}
\sigma _{0}^{\widetilde{\nabla }}(\Delta )=\sigma _{0}^{\nabla }(\Delta
)-\left\langle \sigma _{1}\left( \Delta \right) ,\alpha \right\rangle .
\label{subsymbl}
\end{equation}

\begin{rmk}
The quantizations introduced above are special splitting $R:\mathrm{End}%
\left( \pi \right) \otimes \Sigma _{1}\left( M\right) \rightarrow \mathbf{%
Diff}_{1}\left( \pi ,\pi \right) $~ of \ exact sequence\ (\ref{Exact1}) that
satisfy the following condition:
\begin{equation*}
R\left( A\cdot B\otimes X\right) =A\cdot R\left( B\otimes X\right) ,
\end{equation*}%
for all operators $A,B\in \mathrm{End}\left( \pi \right) $ and vector
fields $X\in \Sigma _{1}\left( M\right) .$ Namely, $\nabla _{X}=R\left(
\mathrm{id}\otimes X\right) .$
\end{rmk}

\section{Pseudogroup actions}

We consider two pseudogroups: $\mathcal{G}\left( M\right) -$ the pseudogroup
of local diffeomorphisms of manifold $M,$ and $\mathbf{Aut}(\xi )-$ the
pseudogroups of local automorphisms of vector bundle $\pi $ over $M.$

There is the following sequence of pseudogroup morphisms
\begin{equation}
1\rightarrow \mathrm{GL}\left( \pi \right) \rightarrow \mathbf{Aut}(\pi
)\rightarrow \mathcal{G}\left( M\right) \rightarrow 1,
\label{group exact seq}
\end{equation}%
where $\mathrm{GL}\left( \pi \right) \subset \mathbf{Aut}(\pi )$ is the
pseudogroup of automorphisms that are identity on $M.$

We'll consider the natural actions of these pseudogroups on sections of the
bundles and operators. Namely, let $\widetilde{\phi }$ be a local
automorphism, $\widetilde{\phi }\in \mathbf{Aut}(\pi ),$ covering a local
diffeomorphism $\phi \in \mathcal{G}\left( M\right) .$

Then we define action of $\widetilde{\phi }$ on sections $s\in C^{\infty
}\left( \xi \right) $ as
\begin{equation*}
\widetilde{\phi }_{\ast }:s\longmapsto \widetilde{\phi }\circ s\circ \phi
^{-1},
\end{equation*}%
and
\begin{equation*}
\widetilde{\phi }_{\ast }:\Delta \longmapsto \widetilde{\phi }_{\ast }\circ
\Delta \circ \widetilde{\phi _{\ast }}^{-1},
\end{equation*}%
for differential operators.

\section{ Artin-Procesi invariants of symbols}

Let's fix a point $a\in M,$ and let $E=\pi _{a}$ be the fibre of the bundle $%
\pi $ at the point $a,$ and let $T=T_{a}\left( M\right) $ be the tangent
space at the point, $\dim T=n\geq 2,\ \dim E=\ m.$ .

In this section we consider symbols of the 1st order differential operators
at the point as tensors of the type $\mathrm{End}\left( E\right) \otimes T$
and classify them with respect to the natural action of the group $G=\mathbf{%
GL}\left( E\right) \times \mathbf{GL}\left( T\right) .$

Remind, that M. Artin made the conjecture (\cite{Art}), concerning
polynomial invariants of $\mathbf{GL}\left( E\right) $-action on $n$-tuples $%
\left( X_{1},...,X_{n}\right) $ of linear operators, $X_{i}\in $ $\mathrm{%
End}\left( E\right) .$ Namely, he stated that the algebra of polynomial
invariants (we denote it by $\mathcal{A}\left( X_{1},...,X_{n}\right) $ ) is
generated by the polynomials $\mathrm{Tr}\left( X_{i_{1}}\cdots
X_{i_{k}}\right) $, associated with non-commutative monomials $\
X_{i_{1}}\cdots X_{i_{k}}.$ This conjecture was proven by C. Procesi in (%
\cite{Pr}) with the important addition that the length of the monomials
could be bounded by $2^{m}-1$ as well as a description of sysygies.

We'll apply this result to our case. Namely, we denote by $\sigma _{\theta
}\in \mathrm{End}\left( E\right) $ the value of symbol $\sigma $ at
covector $\theta \in T^{\ast }$ and consider linear functions $A_{\sigma
,I}:\left( T^{\ast }\right) ^{\otimes k}\rightarrow \mathbb{R}$ , where $%
I=\left( i_{1},...,i_{k}\right) \in S_{k}$ is a permutation of $k$-letters,
as follows%
\begin{equation*}
A_{\sigma ,I}\left( \theta _{1}\otimes \cdots \otimes \theta _{k}\right) =%
\mathrm{Tr}\left( \sigma _{\theta _{i_{1}}}\cdots \sigma _{\theta
_{i_{k}}}\right) .
\end{equation*}%
These tensors $A_{\sigma ,I}\in \left( \left( T^{\ast }\right) ^{\otimes
k}\right) ^{\ast }=T^{\otimes k}$ are obviously $G$-invariants and we call
them \textit{Artin-Procesi tensors}.

Remark also that $A_{\sigma ,I}=A_{\sigma ,J}$, when permutation $J\in S_{k}$
is obtained from $I$ by a cycle permutation.

The following first Artin-Procesi tensors are extremely important for us:

\begin{itemize}
\item In the case $k=1,$ we have a vector $\chi _{\sigma }\in T,$ where $\ $%
\begin{equation*}
\left\langle \theta ,\chi _{\sigma }\right\rangle =\mathrm{Tr}\left( \sigma
_{\theta }\right) ,
\end{equation*}%
for all covectors $\theta \in T^{\ast }.$

\item In the case $k=2,$ we have a symmetric 2-vector $g_{\sigma }\in
S^{2}T, $ where
\begin{equation*}
g_{\sigma }\left( \theta _{1},\theta _{2}\right) =\mathrm{Tr}\left( \sigma
_{\theta _{1}}\sigma _{\theta _{2}}\right) .
\end{equation*}

\item In the case $k=3,$ we have two tensors $h_{i}\in T^{\otimes 3},$ where
\begin{eqnarray*}
h_{1,\sigma }\left( \theta _{1},\theta _{2},\theta _{3}\right) &=&\mathrm{Tr%
}\left( \sigma _{\theta _{1}}\sigma _{\theta _{2}}\sigma _{\theta
_{3}}\right) ,\  \\
h_{2,\sigma }\left( \theta _{1},\theta _{2},\theta _{3}\right) &=&\mathrm{Tr%
}\left( \sigma _{\theta _{2}}\sigma _{\theta _{1}}\sigma _{\theta
_{3}}\right) ,
\end{eqnarray*}%
and equivalently two tensors
\begin{eqnarray*}
h_{\sigma ,s} &=&\frac{1}{2}\left( h_{\sigma ,1}+h_{\sigma ,2}\right) \in
S^{2}T\otimes T,\ \ \  \\
\ \ h_{\sigma ,a} &=&\frac{1}{2}\left( h_{\sigma ,1}-h_{\sigma ,2}\right)
\in \Lambda ^{2}T\otimes T.
\end{eqnarray*}%
All these tensors are $G$-invariants.
\end{itemize}

In what follows we'll use only \textit{symbols in general position}, or
shortly \textit{general symbols}, i.e symbols where the following regularity
conditions hold.

\begin{enumerate}
\item $\ $Vector $\chi _{\sigma }$ \ is non trivial, $\chi _{\sigma }\neq 0.$

\item Quadratic form $g_{\sigma }\in S^{2}T$ on $T^{\ast }$ is non
degenerated.\newline
The inverse quadratic form on $T$ we denote by $g_{\sigma }^{-1}\in
S^{2}T^{\ast },$ and the covector dual to the vector $\chi _{\sigma }$ with
respect to the last quadratic form will be denoted as $\widehat{\chi
_{\sigma }}\in T^{\ast },$ i.e. $\ \widehat{\chi _{\sigma }}\rfloor
g_{\sigma }=\chi _{\sigma },$ or $\chi _{\sigma }\rfloor g_{\sigma }^{-1}=%
\widehat{\chi _{\sigma }}.$

\item Denote also by $\widetilde{h}_{\sigma }\in S^{2}T$ the following
symmetric bivector%
\begin{equation*}
\widetilde{h}_{\sigma }\left( \theta _{1},\theta _{2}\right) =h_{\sigma
,s}\left( \theta _{1},\theta _{2},\widehat{\chi _{\sigma }},\right)
\end{equation*}%
and let $\widehat{h}_{\sigma }:T^{\ast }\rightarrow T^{\ast }$ \ be the
operator, associated with the pair quadratic forms $g_{\sigma }\ $and $%
\widetilde{g}_{\sigma }$ on $T^{\ast }.$ \newline
We'll require that covectors%
\begin{equation*}
e_{1}^{\ast }=\widehat{\chi _{\sigma }},e_{2}^{\ast }=\widehat{h}_{\sigma
}e_{1}^{\ast },...,e_{n}^{\ast }=\widehat{h}_{\sigma }e_{n-1}^{\ast },
\end{equation*}%
are linear independent. \newline
\newline
Then covectors $e^{\ast }\left( \sigma \right) =\left\langle e_{1}^{\ast
},...,e_{n}^{\ast }\right\rangle $ form a coframe in $T^{\ast }$ and vectors
$e\left( \sigma \right) =\left\langle e_{1},...,e_{n}\right\rangle $ of the
dual basis give us the frame in $T.$ Remark that both frames are $G$%
-invariant.

\item In the similar way each covector $\theta \in T^{\ast }\smallsetminus 0$
defines also bivector
\begin{equation*}
\zeta _{\theta ,\sigma }\left( \theta _{1},\theta _{2}\right) =h_{\sigma
,a}\left( \theta _{1},\theta _{2},\theta \right) \in \Lambda ^{2}T.
\end{equation*}%
\newline
Denote by $\Sigma _{\sigma }\subset T^{\ast }$ the variety of  covectors $%
\theta \in T^{\ast }$ where operator $\sigma _{\theta }$ has eigenvalues of
multiplicity greater the one. Then the last condition that $\sigma $ is in
general position requires that $\Sigma _{\sigma }\neq T^{\ast }\ $and $\zeta
_{\theta ,\sigma }$ does not vanish on the complement of $\Sigma _{\sigma }.$
\end{enumerate}

\begin{definition}
A symbol tensor $\sigma \in \mathrm{End}\left( E\right) \otimes T,$ when $%
m\geq 3,$ is said to be in general position if the above conditions hold. In
the case $m=2$ general position requires the first two conditions and that $%
\left\langle \widehat{\chi _{\sigma }},\chi _{\sigma }\right\rangle \neq
0,\Sigma _{\sigma }\neq T^{\ast }.$
\end{definition}

\begin{rmk}
In the case $m=2$ the $G$-invariant frame consist of vectors $e_{1}=\chi
_{\sigma }$ and the vector $e_{2}$ orthonormal to $e_{1}$ with respect to
metric $g_{\sigma }^{-1}.$
\end{rmk}

Let $\sigma $ be a general symbol then, using the invariant frame $e\left(
\sigma \right) ,$ we represent $\sigma $ in the form
\begin{equation*}
\sigma =\sum\limits_{i=1}^{n}\sigma _{i}\otimes e_{i},
\end{equation*}%
where $\sigma _{i}\in \mathrm{End}\left( E\right) .$

Remark that $G$-invariance of the frame \ reduce $G$-equivalence of symbols
to $\mathbf{GL}\left( E\right) $-equivalence of $n$-tuples of operators $%
\left( \sigma _{1},...,\sigma _{n}\right) .$

\begin{proposition}
Let $\sigma $ be a general symbol. Then$\ n\leq m^{2},$ and $\alpha \in
\mathrm{End}\left( E\right) $ is a scalar operator if and only if
commutators $\ [\sigma _{i},\alpha ]=0,$ for all $i=1,...,n$.\newline
\end{proposition}

\begin{pf}
Assume that operators $\sigma _{i}$ are linear dependent, i.e. $%
\sum\limits_{i=1}^{n}\lambda _{i}\sigma _{i}=0,$ for some $\lambda _{i}\in
\mathbb{R\ }$such that $\sum\limits_{i=1}^{n}\lambda _{i}^{2}\neq 0.$ Then $%
\sigma _{\theta }=0,$ for $\theta =\sum\limits_{i=1}^{n}\lambda
_{i}e_{i}^{\ast }\neq 0,$ and $g_{\sigma }\left( \theta ,\theta ^{\prime
}\right) =0,$ for all $\theta ^{\prime }\in T^{\ast }.$ Therefore, $\theta
\in \ker g_{\sigma }$ and $g_{\sigma }$ is degenerated.

Let's now $[\sigma _{i},\alpha ]=0,$ for all $i=1,..,n.$ The 3rd condition
of generality $\sigma $ states that there are covectors $\theta _{1},\theta
_{2}\notin \Sigma _{\sigma },$ such that $[\sigma _{\theta _{1}},\sigma
_{\theta _{2}}]\neq 0.$ Condition $[\sigma _{i},\alpha ]=0,$ for all $i,$
gives $[\sigma _{\theta _{1}},\alpha ]=[\sigma _{\theta _{2}},\alpha ]=0.$
Therefore operator $\alpha $ has the same eigenvectors with operators $%
\sigma _{\theta _{1}}\ $and $\ \sigma _{\theta _{2}},$ that possible only
when $\alpha $ is a scalar operator.
\end{pf}

\begin{corollary}
The stationary algebra Lie of a general symbol consist of scalar operators
and codimension of $G$-orbit of such symbols equals to
\begin{equation*}
\nu =\left( n-1\right) \left( m^{2}-n-1\right) .
\end{equation*}
\end{corollary}

Take now Artin-Procesi tensors $A_{\sigma ,I}\in T^{\otimes k},k=\left\vert
I\right\vert ,$ for general symbol $\sigma $, and write down them in the
invariant frame $e\left( \sigma \right) :$%
\begin{equation*}
A_{\sigma ,I}=\sum\limits_{J}A_{I}^{J}\left( \sigma \right) e^{J},
\end{equation*}%
where $J=\left( j_{1},...,j_{k}\right) $ and $e^{J}=e_{j_{1}}\otimes \cdots
\otimes e_{j_{k}}.$

Then coefficients $A_{I}^{J}\left( \sigma \right) $ are rational functions
on $\mathrm{End}\left( E\right) \otimes T$ and $G$-invariants. We call them
\textit{Artin-Procesi invariants} of the symbol.

We have the algebraic $G$-action on $\mathrm{End}\left( E\right) \otimes T.$
Therefore, due to Rosenlicht theorem (\cite{Ros}), there are $\nu $
Artin-Procesi invariants $a_{1},..,a_{\nu }$ which are algebraically
independent and all rational $G$-invariants (in particular all Artin-Procesi
invariants) are rational functions of them.

In this case we call $a_{1},..,a_{\nu }$ \textit{basic Artin-Procesi
invariants}.

$G$-orbits $G\left( \sigma \right) \subset $ $\left\{ a_{1}=\mathrm{const}%
_{1},..,a_{\nu }=\mathrm{const}_{\nu }\right\} $ such that differentials $%
da_{1},..,da_{\nu }$ are linear independent at the points of \ the orbit we
call \textit{regular} as well as the symbol $\sigma $ itself.

Remark that $G$-orbits are connected and therefore the manifold $\left\{
a_{1}=\mathrm{const}_{1},..,a_{\nu }=\mathrm{const}_{\nu }\right\} $ in
the last case is a finite union of regular orbits.

\begin{example}
The following invariants might be used as basic:

\begin{enumerate}
\item For $n=2$
\begin{eqnarray*}
I_{i}^{\left( 1\right) } &=&\mathrm{Tr}\sigma _{1}^{i},I_{i}^{\left(
2\right) }=\mathrm{Tr}\sigma _{2}^{i},i=1,...,m, \\
I_{i,j} &=&\mathrm{Tr}\left( \sigma _{1}^{i}\sigma _{2}^{j}\right) ,\ 1\leq
i,j\leq m-1.
\end{eqnarray*}

\item For $n\geq 3$
\begin{eqnarray*}
I_{i}^{\left( k\right) } &=&\mathrm{Tr}\sigma _{k}^{i},\ \
k=1,..n,i=1,...,m, \\
I_{i,j}^{\left( k\right) } &=&\mathrm{Tr}\left( \sigma _{1}^{i}\sigma
_{k}^{j}\right) ,\ \ \ 1\leq i,j\leq m-1, \\
J_{j}^{\left( k\right) } &=&\mathrm{Tr}\left( \sigma _{1}\sigma _{2}\sigma
_{k}^{j}\right) ,\ \ k=3,..n,j=1,...,m-1.
\end{eqnarray*}
\end{enumerate}
\end{example}

Summarizing we get the following.

Let $\sigma \in \mathrm{End}\left( E\right) \otimes T$ be a general symbol.
\ Then

\begin{enumerate}
\item Tensors
\begin{equation*}
A_{\sigma ,I}=\mathrm{Tr}\left( \sigma _{\theta _{i_{1}}}\cdots \sigma
_{\theta _{i_{k}}}\right) \in T^{\otimes \left\vert I\right\vert },
\end{equation*}%
are $\mathbf{GL}\left( E\right) \times \mathbf{GL}\left( T\right) $%
-invariants.

\item The conditions that $\sigma $ is in general position are conditions on
tensors $A_{\sigma ,I},$ with $\left\vert I\right\vert \leq 3.$ Under these
conditions the symbol $\sigma $ determines $\mathbf{GL}\left( E\right)
\times \mathbf{GL}\left( T\right) $-invariant coframe $e\left( \sigma
\right) $ in the tangent space $T.$

\item Coefficients  $A_{\sigma ,I}^{J}$ of the tensors $A_{\sigma ,I}$ in
the frame $e\left( \sigma \right) $ are $\mathbf{GL}\left( E\right) \times
\mathbf{GL}\left( T\right) $-invariants of symbols. These invariants are
polynomials over the field $\mathbb{R}\left( \det g_{\sigma }^{-1}\right) .$

\item The maximum dimensional $\mathbf{GL}\left( E\right) \times \mathbf{GL}%
\left( T\right) $-orbits in $\mathrm{End}\left( E\right) \otimes T$ have
codimension $\nu =\left( n-1\right) \left( m^{2}-n-1\right) .$ The $\mathbf{%
GL}\left( E\right) \times \mathbf{GL}\left( T\right) $-action is algebraic
and, due to the Rosenlicht theorem (\cite{Ros}) the field of rational $%
\mathbf{GL}\left( E\right) \times \mathbf{GL}\left( T\right) $-invariants
has transcendent degree $\nu .$ We call $\mathbf{GL}\left( E\right) \times
\mathbf{GL}\left( T\right) $-orbit (and its elements) \textit{regular} if
there are \textit{basic} invariants $a_{1},..,a_{\nu }$ such that their
differentials are linear independent at the points of the orbit.
\end{enumerate}

\section{Connections, associated with regular symbols}

Let $\sigma =\sigma _{1}\left( \Delta \right) $ be a regular symbol and let $%
\nabla _{c}$ be the Levi-Civita connection in the tangent bundle $\tau
\left( M\right) ,$ associated with the metric $g=g_{\sigma }.$

Let $\nabla $ be a connection in the bundle $\pi .$ We'll use the same
notation $\nabla $ for the associated connection in the bundle $\mathrm{End}%
\left( \pi \right) ,$ where $\nabla _{X}\left( A\right) =\nabla _{X}\circ
A-A\circ \nabla _{X},$ for all $X\in \Sigma _{1}\left( M\right) ,\ A\in
\mathrm{End}\left( \pi \right) .$ \

Let $\nabla ^{c}=\nabla \otimes \nabla _{c},$ be the connection in the
bundle $\mathrm{End}\left( \pi \right) \otimes \Sigma _{1}\left( M\right) $
equals to the tensor product of the connection $\nabla $ in the bundle $%
\mathrm{End}\left( \pi \right) $ and the Levi-Civita connection in the
tangent bundle $\tau (M).$

The covariant differential $d_{\nabla ^{c}}:\mathrm{End}\left( \pi \right)
\otimes \Sigma _{1}\left( M\right) \rightarrow \mathrm{End}\left( \pi
\right) \otimes \Sigma _{1}\left( M\right) \otimes \Omega ^{1}\left(
M\right) =\mathrm{End}\left( \pi \right) \otimes \mathrm{End}\left( \tau
\right) ,$ applying to the symbol $\sigma ,$ gives us a tensor
\begin{equation*}
d_{\nabla ^{c}}\sigma \in \mathrm{End}\left( \pi \right) \otimes \mathrm{%
End}\left( \tau \right) .
\end{equation*}%
Remark that for another connection $\widetilde{\nabla }$ in the bundle $\pi $
such that
\begin{equation*}
\widetilde{\nabla }_{X}-\nabla _{X}=\alpha \left( X\right) \in \mathrm{End}%
\left( \pi \right) ,
\end{equation*}%
where $\alpha \in \mathrm{End}\left( \pi \right) \otimes \Omega ^{1}\left(
M\right) ,$ we have
\begin{equation*}
d_{\widetilde{\nabla }^{c}}\sigma -d_{\nabla ^{c}}\sigma =[\alpha ,\sigma ],
\end{equation*}%
where by $[\alpha ,\sigma ]\in \mathrm{End}\left( \pi \right) \otimes
\mathrm{End}\left( \tau \right) $ we denoted the natural pairing%
\begin{eqnarray*}
\lbrack ,]:\mathrm{End}\left( \pi \right) \otimes \Omega ^{1}\left(
M\right) \times \mathrm{End}\left( \pi \right) \otimes \Sigma _{1}\left(
M\right) \rightarrow \mathrm{End}\left( \pi \right) \otimes \mathrm{End}%
\left( \tau \right) &&, \\
\lbrack A\otimes \theta ,B\otimes X]=[A,B]\otimes \left( \theta \otimes
X\right) \in \mathrm{End}\left( \pi \right) \otimes \mathrm{End}\left(
\tau \right) &&.
\end{eqnarray*}%
Indeed, let $\sigma =\Sigma _{i}A_{i}\otimes X_{i},$ where $A_{i}\in
\mathrm{End}\left( \pi \right) ,$ and $X_{i}\in \Sigma _{1}\left( M\right)
. $ Then%
\begin{equation*}
\nabla _{X}^{c}\left( \sigma \right) =\Sigma _{i}\left( \nabla _{X}\left(
A_{i}\right) \otimes X_{i}+A_{i}\otimes \nabla _{cX}\left( X_{i}\right)
\right)
\end{equation*}%
and
\begin{eqnarray*}
\widetilde{\nabla }_{X}^{c}\left( \sigma \right) =\Sigma _{i}\left(
\widetilde{\nabla }_{X}^{c}\left( A_{i}\right) \otimes X_{i}+A_{i}\otimes
\nabla _{cX}\left( X_{i}\right) \right) = && \\
=\nabla _{X}^{c}\left( \sigma \right) +\Sigma _{i}[\alpha \left( X\right)
,A_{i}]\otimes X_{i}=\nabla _{X}^{c}\left( \sigma \right) +[\alpha ,\sigma
]\left( X\right) . &&
\end{eqnarray*}

\begin{proposition}
Operator $D\left( \sigma \right) =\left( \mathrm{Tr}\otimes \mathrm{id}\right)
\left( d_{\nabla ^{c}}\sigma \right) \in \mathrm{End}\left( \tau \right) $
does not depend on choice of \ connection $\nabla .$
\end{proposition}

Let now $\sigma $ be a general symbol and let $\sigma =\sum_{i}\sigma
_{i}\otimes e_{i}$ be its decomposition in the invariant frame. Remark that $%
d_{\widetilde{\nabla }^{c}}\sigma =d_{\nabla ^{c}}\sigma $ if and only if $%
[\alpha ,\sigma ]=0$ or if $[\alpha _{i},\sigma _{j}]=0,$ for all $i,j,$
where $\alpha =\sum_{j}\alpha _{i}\otimes e_{i}^{\ast }.$ Therefore, each
operator $\alpha _{i}$ is a scalar operator $\alpha _{i}=\lambda _{i}$ $%
\mathrm{id}.$ In other words,
\begin{equation*}
d_{\widetilde{\nabla }^{c}}\sigma =d_{\nabla ^{c}}\sigma \ \Leftrightarrow
\widetilde{\nabla }-\nabla =\mathrm{id}\otimes \lambda ,
\end{equation*}%
where $\lambda \in \Omega ^{1}\left( M\right) $ is a differential $1-$form
such that $\left\langle \lambda ,e_{i}\right\rangle =\lambda _{i}.$

Denote by $W\left( \sigma \right) \subset \mathrm{End}\left( \pi \right)
\otimes \mathrm{End}\left( \tau \right) $ the affine subbundle generated by
tensors $d_{\nabla ^{c}}\sigma ,$ computing for all connections $\nabla
^{c}. $ The fibres of this bundle isomorphic to $\mathrm{sl}\left( \pi
\right) \otimes \mathrm{End}\left( \tau \right) ,$ where%
\begin{equation*}
\mathrm{sl}\left( \pi \right) =\left\{ \left. A\in \mathrm{End}\left( \pi
\right) \right\vert \ \mathrm{Tr}A=0\right\} .
\end{equation*}

Let's consider $\mathrm{End}\left( \pi \right) \otimes \mathrm{End}\left(
\tau \right) $ as Euclidean bundle with respect to inner product $\left(
A\otimes a,B\otimes b\right) =\mathrm{Tr}\left( AB\right) \cdot \mathrm{Tr}%
\left( ab\right) ,$ where $A,B\in \mathrm{End}\left( \pi \right) $ and $%
a,b\in \mathrm{End}\left( \tau \right) .$

A connection $\nabla $ we call \textit{minimal }if tensor $d_{\nabla
^{c}}\sigma $ orthogonal to $W\left( \sigma \right) $ with respect to the
above inner product.

\begin{proposition}
Let $\sigma $ be a symbol in general position. Then, for any connection $%
\nabla $ in the bundle $\pi ,$ there is and unique tensor $\alpha \in
\mathrm{sl}\left( \pi \right) \otimes \Omega ^{1}\left( M\right) ,$ such
that
\begin{equation}
\left( d_{\nabla ^{c}}\sigma +[\alpha ,\sigma ],\mathrm{sl}\left( \pi
\right) \otimes \mathrm{End}\left( \tau \right) \right) =0.  \label{Alpha}
\end{equation}
\end{proposition}

\begin{pf}
Let $E_{ij}=e_{i}\otimes e_{j}^{\ast }$ be the elementary operators in the
invariant basis in $T$ and let%
\begin{equation*}
d_{\nabla ^{c}}\sigma -\mathrm{id}\otimes D\left( \sigma \right)
=\sum_{i,j}d_{ij}\otimes E_{ij},
\end{equation*}%
for some $d_{ij}\in \mathrm{sl}\left( \pi \right) ,$ and%
\begin{equation*}
\sigma =\sum_{i}\sigma _{i}\otimes e_{i},\ \ \alpha =\sum_{i}\alpha
_{j}\otimes e_{j}^{\ast },
\end{equation*}%
where $\alpha _{j}\in \mathrm{sl}\left( \pi \right) .$

Then (\ref{Alpha}) equivalent to the following linear system of $n\left(
m^{2}-1\right) $ equations%
\begin{equation}
d_{ij}-[\alpha _{j},\sigma _{i}]=0,  \label{Alpha2}
\end{equation}%
with respect to $n\left( m^{2}-1\right) $ unknowns $\alpha _{j}.$

As we have seen above the stationary Lie algebra of the symbol consists of
scalar operators. It means that the homogeneous system of (\ref{Alpha2}) has
the only zero solution.
\end{pf}

The above theorem states the existence of \ minimal connections and, as we
have seen, any two minimal connections differ on tensors of the form $\alpha
=\mathrm{id}\otimes \lambda .$

Let $R_{\nabla }\left( \sigma \right) \in \mathrm{End}\left( \pi \right)
\otimes \Omega ^{2}\left( M\right) $ be the curvature tensor of minimal
connection $\nabla .$ Then (see, for example, \cite{Hus}) we have
\begin{equation*}
R_{\widetilde{\nabla }}\left( \sigma \right) -R_{\nabla }\left( \sigma
\right) =\mathrm{id}\otimes d\lambda ,
\end{equation*}%
if $\ \widetilde{\nabla }-\nabla =\mathrm{id}\otimes \lambda .$

Therefore, the tensor
\begin{equation*}
R^{0}\left( \sigma \right) =R_{\nabla }\left( \sigma \right) -\frac{1}{m}%
\mathrm{id}\otimes \mathrm{ch}\left( \nabla \right) ,
\end{equation*}%
where $\mathrm{ch}\left( \nabla \right) =\mathrm{Tr}\left( R_{\nabla
}\left( \sigma \right) \right) \in \Omega ^{2}\left( M\right) $ is the first
\textit{Chern form}, does not depend on choice of minimal connection $\nabla
.$

This tensor is also $G$-invariant of the symbol$.$

On the other hand,we have (\ref{subsymbl})
\begin{equation*}
\sigma _{0}\left( \widetilde{\nabla }\right) -\sigma _{0}\left( \nabla
\right) =\sigma _{-\lambda },
\end{equation*}%
for given operator $\Delta \in \mathbf{Diff}_{1}\left( \pi ,\pi \right) $
and minimal connections $\nabla $ and $\widetilde{\nabla }\ .$

Therefore, subbundle $W_{0}\left( \Delta \right) \subset \mathrm{End}\left(
\pi \right) ,$ generated by operators $\sigma _{0}\left( \nabla \right) ,$
taking for all minimal connections $\nabla $, is an affine bundle isomorphic
to $\tau ^{\ast }$ by the injection $\lambda \in \Omega ^{1}\rightarrow
\sigma _{\lambda }\in \mathrm{End}\left( \pi \right) .$

We say that a minimal connection $\nabla $ is \textit{associated with
operator} $\Delta $ if the operator $\sigma _{0}\left( \nabla \right) $ is
orthogonal to $W_{0}\left( \Delta \right) .$

To find this connection we should find differential form $\lambda $ such
that
\begin{equation*}
\mathrm{Tr}\left( \left( \sigma _{0}\left( \nabla \right) -\sigma _{\lambda
}\right) \sigma _{\theta }\right) =0,
\end{equation*}%
for all $\theta \in \Omega ^{1}\left( M\right) .$

Remark that the last system has a unique solution because quadratic form $%
g_{\sigma }\left( \lambda ,\theta \right) =\mathrm{Tr}\left( \sigma _{\lambda
}\sigma _{\theta }\right) $ is non degenerated.

Finally we get the following result.

\begin{thm}
Let $\Delta \in \mathbf{Diff}_{1}\left( \pi ,\pi \right) $ be a differential
operator with general symbol $\sigma \in \mathrm{End}\left( \pi \right)
\otimes \Sigma _{1}$.. Then there exists and unique associated connection $%
\nabla ^{\Delta }$ in the bundle $\pi $ which is minimal and the operator $%
\sigma _{0}\left( \nabla ^{\Delta }\right) \in \mathrm{End}\left( \pi
\right) $ in the decomposition
\begin{equation*}
\Delta =Q_{\nabla ^{\Delta }}\left( \sigma \right) +\sigma _{0}\left( \nabla
^{\Delta }\right)
\end{equation*}%
satisfies the following conditions:%
\begin{equation*}
\mathrm{Tr}\left( \sigma _{0}\left( \nabla ^{\Delta }\right) \cdot \sigma
_{\theta }\right) =0,
\end{equation*}%
for all differential 1-forms $\theta \in \Omega ^{1}\left( M\right) .$
\end{thm}

\begin{corollary}
Let $A\in \mathrm{Aut}\left( \pi \right) $ be an automorphism, transforming
operator $\Delta \in \mathbf{Diff}_{1}\left( \pi ,\pi \right) $ to operator $%
\Delta ^{\prime }\in \mathbf{Diff}_{1}\left( \pi ,\pi \right) ,$ $A_{\ast
}\left( \Delta \right) =\Delta ^{\prime },$ and let $\Delta =Q_{\nabla
^{\Delta }}\left( \sigma \right) +\sigma _{0}\left( \nabla ^{\Delta }\right)
$ and $\Delta ^{\prime }=Q_{\nabla ^{\Delta ^{\prime }}}\left( \sigma
^{\prime }\right) +\sigma _{0}\left( \nabla ^{\Delta \prime }\right) $ be
the decompositions with respect to the associated connections. \newline
Then $A_{\ast }\left( \sigma \right) =\sigma ^{\prime },A_{\ast }\left(
\nabla ^{\Delta }\right) =\nabla ^{\Delta ^{\prime }},$ and \ $A_{\ast
}\left( \sigma _{0}\left( \nabla ^{\Delta }\right) \right) =\sigma
_{0}\left( \nabla ^{\Delta \prime }\right) .$
\end{corollary}

\begin{pf}
The statement follows directly from the fact that any automorphism $A\in
\mathrm{GL}\left( \pi \right) $ preserves the inner product structure in $%
\mathrm{End}\left( \pi \right) .$
\end{pf}

This corollary shows that the decomposition with respect to associated
connection behave in the natural way under transformations from the
pseudogroup $\mathrm{Aut}\left( \pi \right) .$ It is allow us to extend
Artin invariants of symbols to \textit{Artin-Procesi invariants of operators.%
}

Namely, let $\Delta \in \mathbf{Diff}_{1}\left( \pi ,\pi \right) $ be a
differential operator with regular symbol at a point $a\in M.$ Then the
decomposition of the symbol
\begin{equation*}
\sigma =\Sigma _{i}\sigma _{i}\otimes e_{i}
\end{equation*}
in the invariant frame $e\left( \sigma \right) $ and operator $\sigma _{0}$
gives us the Artin-Procesi invariants of operators as elements of algebra $%
\mathcal{A}\left( \sigma _{0},\sigma _{1}\right) .$

Moreover, the curvature tensor
\begin{equation*}
R_{\nabla ^{\Delta }}\left( \sigma \right) \in \mathrm{End}\left( \pi
\right) \otimes \Omega ^{2}\left( M\right) ,
\end{equation*}%
computing for the associated connection $\nabla ^{\Delta }$ in the bundle $%
\pi $ is an $\mathbf{Aut}(\pi )$-invariant of the operators as well as the
closed Chern differential 2-form
\begin{equation*}
\mathrm{ch}\left( \Delta \right) =\mathrm{Tr}\left( R_{\nabla ^{\Delta
}}\left( \sigma \right) \right) \in \Omega ^{2}\left( M\right) ,
\end{equation*}%
which we call \textit{Chern invariant of the operator}.

Coefficients of the form in the invariant frame
\begin{equation*}
\mathrm{ch}_{ij}\left( \Delta \right) =\mathrm{ch}\left( \Delta
\right) \left( e_{i},e_{j}\right)
\end{equation*}%
are scalar $\mathbf{Aut}(\pi )$-invariants of the the operators.

\section{Equivalence of regular differential operators}

A differential operator $\Delta \in \mathbf{Diff}_{1}\left( \pi ,\pi \right)
$ is said to be \textit{regular} at a point $a\in M$ \ if its symbol $\sigma
\in \mathrm{End}\left( \pi \right) \otimes \Sigma _{1}\left( M\right) $ is
regular\ at the point and among Artin-Procesi invariants $a_{1},...,a_{\nu
_{0}}$ \ of the operator, defining the $\mathrm{GL}\left( \pi \right)
\times \mathrm{GL}\left( T\right) \ -$orbit of the pairs $\left( \sigma
_{0},\sigma \right) ,$ there are $n=\dim M$ \ invariants, say $%
a_{1},..,a_{n},$ such that differentials of functions $a_{1}\left( \Delta
\right) ,..,a_{n}\left( \Delta \right) $ are linear independent at the point
$a\in M.$

In this case there is a neighborhood $a\in U\subset M$ such that $\mathrm{GL%
}\left( \pi _{b}\right) \times \mathrm{GL}\left( T_{b}\right) \ -$orbits of
$\ $the pairs $\left( \sigma _{0}\left( b\right) ,\sigma \left( b\right)
\right) \in \mathrm{End}\left( \pi _{b}\right) \oplus $ $\mathrm{End}%
\left( \pi _{b}\right) \otimes T_{b}$\ are regular and are defined by values
the same basic invariants $a_{1},...,a_{\nu _{0}},$ and furthermore
functions $a_{1}\left( \Delta \right) ,..,a_{n}\left( \Delta \right) $ are
local coordinates in $U.$

We call such local coordinates \textit{natural}.

\begin{rmk}
The regularity of symbol requires that $m^{2}\geq n+1,$ and the number $\nu
_{0}$ of basic Artin-Procesi invariants of operators equals%
\begin{equation*}
\nu _{0}=\nu +m^{2}>n.
\end{equation*}
\end{rmk}

Let's $a_{1},...,a_{n}$ be the natural coordinates then all basic
Artin-Procesi invariants $a_{j}\left( \Delta \right) $ , $1\leq j\leq \nu
_{0},$ are functions of $a_{1}\left( \Delta \right) ,..,a_{n}\left( \Delta
\right) $ in $U:$
\begin{equation}
a_{j}\left( \Delta \right) =F_{j}\left( a_{1}\left( \Delta \right)
,..,a_{n}\left( \Delta \right) \right) ,\ \ n+1\leq j\leq \nu _{0},
\label{loc}
\end{equation}%
for some functions $F_{j}.$

All rational $G$-invariants of the pairs $\left( \sigma _{0},\sigma \right) $
are rational functions of basic Artin-Procesi invariants $\left(
a_{1},..,a_{\nu _{0}}\right) $ and, due to (\ref{loc}), their values in the
neighborhood $U$ are also functions of $\left( a_{1}\left( \Delta \right)
,..,a_{n}\left( \Delta \right) \right) $ completely defined by functions $%
F_{j}.$

Let's$\ \phi _{U}:U\rightarrow \mathbf{D}_{U}\subset \mathbb{R}^{n}$ be the
natural local chart,
\begin{equation*}
\phi _{U}\left( b\right) =\left( a_{1}\left( \Delta \right) \left( b\right)
,..,a_{n}\left( \Delta \right) \left( b\right) \right) ,
\end{equation*}
and let $F_{U}:\mathbf{D}_{U}\rightarrow \mathbb{R}^{\nu _{0}-n}$ be the
function given by (\ref{loc}), and
\begin{equation*}
\mathrm{ch}_{U}\left( \Delta \right) =\phi _{U\ast }\left(
\mathrm{ch}\left( \Delta \right) \right) \in \Omega ^{2}\left( \mathbf{D}%
_{U}\right) .
\end{equation*}

The data $\left( \phi _{U},\mathbf{D}_{U},F_{U},\mathrm{ch}_{U}\left(
\Delta \right) \right) $ we call \textit{model of the differential operator}
in neighborhood $U.$

\begin{thm}
Let differential opertors $\Delta $ and $\Delta ^{\prime }$ has the same
model in a simply connected open set $U.$ \newline
Then there is and a unique automorphism $A_{U}\in \mathrm{GL}\left( \pi
_{U}\right) $ of the restriction bundle $\pi $ on\ domain $U$ such that $%
A_{U\ast }\left( \Delta \right) =\Delta ^{\prime }.$
\end{thm}

\begin{pf}
Condition that two operators have the same model means that the pairs $%
\left( \sigma _{0},\sigma \right) $ and $\left( \sigma _{0}^{\prime },\sigma
^{\prime }\right) $ belong to the same $G$-orbit at any point $a\in U.$
Therefore, there are the above automorphisms $A_{U}.$ Moreover, if $\nabla $
is the connection, associated with operator $\Delta ,$ then $\nabla ^{\prime
}=A_{U\ast }(\nabla )$ is a minimal connection for $\Delta ^{\prime }$
because $A_{U\ast }$ preserves the inner structure in $\mathrm{End}\left(
\pi \right) \otimes \mathrm{End}(T)$ and transforms $\sigma $ to $\sigma
^{\prime }.$

Let
\begin{equation*}
\nabla ^{\Delta ^{\prime }}=\nabla ^{\prime }+\mathrm{id}\otimes \lambda
\end{equation*}%
for some differential 1-form $\lambda .$

Then, $\mathrm{ch}\left( \Delta ^{\prime }\right) =\mathrm{ch}\left(
\nabla ^{\prime }\right) +d\lambda =\mathrm{ch}\left( \Delta \right)
+d\lambda ,$ and $d\lambda =0$. Therefore, $\lambda =df$ in $U.$

Finally, $\left( B_{U}\right) \left( \nabla \right) =\nabla ^{\Delta
^{\prime }},\ $where $B_{U}=$ $e^{-f}A_{U},$ and therefore $B_{U\ast }\left(
\Delta \right) =\Delta ^{\prime }.$
\end{pf}

\begin{corollary}
Let differential opertors $\Delta $ and $\Delta ^{\prime }$ has models $%
\left( \phi _{U},\mathbf{D}_{U},F_{U},\mathrm{ch}_{U}\left( \Delta \right)
\right) $ and $\left( \phi _{U^{\prime }},\mathbf{D}_{U^{\prime
}},F_{U^{\prime }},\mathrm{ch}_{U}\left( \Delta ^{\prime }\right) \right) $
in open sets $U\subset M$ and $U^{\prime }\subset M$ \ defined by the same
basic invariants $a_{1},..,a_{n}$ and functions $F_{j}.$ \newline
Let $\widetilde{U}\subset U$ and $\widetilde{U^{\prime }}\subset U^{\prime }$
be open and simply connected domains such that $\phi _{U}\left( \widetilde{U}%
\right) =\phi _{U^{\prime }}\left( \widetilde{U^{\prime }}\right) \subset
\mathbf{D}_{U}\cap \mathbf{D}_{U^{\prime }},$ and where $\mathrm{ch}%
_{U}\left( \Delta \right) =\mathrm{ch}_{U^{\prime }}\left( \Delta ^{\prime
}\right) .$ \newline
Then there is an automorphism $A_{\widetilde{U},\widetilde{U^{\prime }}}:\pi
_{\widetilde{U}}\rightarrow \pi _{\widetilde{U^{\prime }}}$ \ such that $%
\left( A_{\widetilde{U},\widetilde{U^{\prime }}}\right) _{\ast }\Delta
=\Delta ^{\prime }.$
\end{corollary}

Let now $\Delta $ be a differential operator regular on the manifold $M,$
i.e.regular at all points of the manifold. We say that an atlas $\left\{
\left( \phi _{U^{\alpha }},\mathbf{D}_{U^{\alpha }}\right) \right\} $ given
by models $\left\{ \left( \phi _{U^{\alpha }},\mathbf{D}_{U^{\alpha
}},F_{U^{\alpha }},\mathrm{ch}_{U^{\alpha }}\left( \Delta \right) \right)
\right\} $ is \textit{natural} \ if the\ sets of basic invariants $\left(
a_{1}^{\alpha },...,a_{n}^{\alpha }\right) $ are different for different $%
\alpha .$

The following result follows directly from the above theorem.

\begin{thm}
Two linear differential operators $\Delta ,\Delta ^{\prime }\in \mathbf{Diff}%
_{1}\left( \pi ,\pi \right) $ on a manifold $M$ are $\mathrm{Aut}\left( \pi
\right) $-equivalent if and only if a natural atlas for  operator $\Delta $
is the natural atlas for $\Delta ^{\prime },$ i.e. they have the same models.
\end{thm}

\section*{Acknowledgements}
The author was partially supported by the Russian Foundation for Basic
Research (project 18-29-10013).



\end{document}